# ASYMPTOTIC THEORY OF LEAST SQUARES ESTIMATORS FOR NEARLY UNSTABLE PROCESSES UNDER STRONG DEPENDENCE

By Boris Buchmann and Ngai Hang Chan

*Australian National University and Chinese University of Hong Kong*

This paper considers the effect of least squares procedures for nearly unstable linear time series with strongly dependent innovations. Under a general framework and appropriate scaling, it is shown that ordinary least squares procedures converge to functionals of fractional Ornstein–Uhlenbeck processes. We use fractional integrated noise as an example to illustrate the important ideas. In this case, the functionals bear only formal analogy to those in the classical framework with uncorrelated innovations, with Wiener processes being replaced by fractional Brownian motions. It is also shown that limit theorems for the functionals involve nonstandard scaling and nonstandard limiting distributions. Results of this paper shed light on the asymptotic behavior of nearly unstable long-memory processes.

**1. Introduction.** Consider a first-order autoregressive model,

$$X_t = \beta X_{t-1} + \varepsilon_t \qquad \text{for } t \in \mathbb{N}, X_0 = 0.$$

The parameter $\beta$ is unknown and has to be estimated from the observations $X_1, \ldots, X_n$. Whenever $\varepsilon = (\varepsilon_t)_n$ is a sequence of independent standard normal random variables independent of $X_0$ both the least squares estimator and the maximum likelihood estimator for $\beta_n$ are given by the formula

$$\hat{b}_n = \hat{b}_n(X_0, \ldots, X_n) = \frac{\sum_{t=0}^{n-1} X_{t+1} X_t}{\sum_{t=0}^{n-1} X_t^2}.$$

For $|\beta| < 1$ (stationary regime), Mann and Wald [18] showed that

$$(1.1) \qquad \hat{\tau}_n = \left(\sum_{t=0}^{n-1} X_t^2\right)^{1/2} (\hat{b}_n - \beta) \stackrel{d}{\to} W_1, \qquad n \to \infty,$$

where $W_1$ is a standard normal random variable and $\stackrel{d}{\to}$ denotes convergence in distribution as $n \to \infty$. For the explosive case $|\beta| > 1$, (1.1) holds only when $\varepsilon$ is a sequence of independent standard normal random variables; see Anderson [1]. In general, if $\varepsilon$ does not form a sequence of identically distributed and independent random variables, a limit distribution for $\hat{\tau}_n$ may not even exist. On the other hand, (1.1) fails to hold for $\beta = 1$ even when $\varepsilon$ forms a sequence of independent standard normal random variables. Surprisingly, the limit distribution of $\hat{\tau}_n$ does exist, but is not normal. White [30] and Rao [21] showed that it is a functional of Brownian motion $W = (W_t)_{0 \le t \le 1}$, that is,

$$(1.2) \qquad \hat{\tau}_n \stackrel{d}{\to} \bar{\tau} = \tfrac{1}{2}[W_1^2 - 1]\left[\int_0^1 W_s^2 \, ds\right]^{-1/2}, \qquad n \to \infty.$$

This contrast between (1.1) and (1.2) is an example of a critical phenomenon. The parameter value $\beta = 1$ comprises a singularity and there is a lack of smooth transition of the limiting distribution of $\hat{\tau}_n$ when $\beta$ is close to one. In particular, for finite sample analysis or tests under local alternatives, a key question becomes that if $\beta$ is close to one, what kind of approximation should be used for $\hat{\tau}_n$? An answer to this question is given in Chan and Wei [4], where a class of nearly nonstationary models is proposed. In the econometric literature, such a class is also known as the near-integrated time series; see Phillips [20]. In a spectral setting, Dahlhaus [5] consider tapered and nontapered Yule–Walker estimates of near-integrated models. Let $n \in \mathbb{N}$ and consider

$$(1.3) \qquad X_t^{(n)} = \beta_n X_{t-1}^{(n)} + \varepsilon_t \qquad \text{for } t = 1, \ldots, n,$$

where $X_0^{(n)}$ is an initial value.

THEOREM 1.1. (i) *Let $X_0^{(n)} = 0$, almost surely, for all $n \in \mathbb{N}$. Suppose that $\varepsilon$ forms a sequence of martingale differences such that for $n \to \infty$*

$$(1.4) \qquad \frac{1}{n}\sum_{t=2}^{n} E[\varepsilon_t^2 | \mathcal{F}_{t-1}] = 1 + o_P(1),$$

$$(1.5) \qquad \frac{1}{n}\sum_{t=2}^{n} E[\varepsilon_t^2 1_{|\varepsilon_t| > n^{1/2}\alpha} | \mathcal{F}_{t-1}] = o_P(1) \qquad \text{for all } \alpha > 0,$$

*where $\mathcal{F}_t = \sigma(\varepsilon_s : 0 \le s \le t)$.*

*If there exists $\gamma \in \mathbb{R}$ such that $\beta_n = 1 - \gamma/n$, then*

$$\hat{\tau}_n \xrightarrow{d} \bar{\tau}(\gamma) = \frac{\int_0^1 W_{\gamma,s}\, dW_s}{\sqrt{\int_0^1 W_{\gamma,s}^2\, ds}}, \qquad n \to \infty,$$

*where $W_{\gamma,t} = W_t - \gamma \int_0^t e^{-\gamma(t-s)} W_s\, ds$, $t \in [0,1]$.*

(ii) $\bar{\tau}_0 \stackrel{d}{=} \bar{\tau}$ *and* $\mathcal{L}(\bar{\tau}_\gamma) \xrightarrow{d} \mathcal{L}(W_1)$ *for $|\gamma| \to \infty$.*

Model (1.3) was later explored by Jeganathan [14, 15], who generalized the near-integrated notion to a general AR($p$) case. In practice, econometric and financial time series often exhibit long-range dependent structure (see, e.g., Robinson [22, 24] and Doukhan, Oppenheimer and Taqqu [11]) which cannot be encompassed by the martingale difference setting of Chan and Wei [4].

In this paper we are concerned with least squares estimators when long-range dependence in the innovation is present. We investigate the behavior of the transition of the limiting distributions of the least squares statistic $\hat{\tau}_n$ when $\beta$ is close to one, something similar to (1.1) and (1.2). We apply our results to noise sequences $\varepsilon$ with corresponding partial sum processes lying in the domain of attraction of fractional Brownian motions. Under this setting, we investigate the behavior of the transition of the limiting distributions. Although a formal analogy will be developed, as will be seen, the techniques and the proofs involved are very different between the long-range and short-range dependent cases. In particular, our theory relies on making use of the fractional Ornstein–Uhlenbeck processes.

It should be pointed out that there exist other methods for dealing with long memory models, notably the Whittle estimates, the tapered Whittle estimates and related semiparametric procedures that were discussed in Velasco and Robinson [29] and Robinson [24]. For further details and properties about these methods, we refer the readers to these papers and the references therein. The semiparametric approach constitutes an important alternative to ordinary least squares (OLS) as its origin stems from a basic whitening idea, which could be argued to be one of the major goals in modern time series analysis. Another important advantage of the semiparametric procedure is that their estimators are usually asymptotically normal and thus standard limit theory can be applied so that optimality properties and relative neatness can be attained, albeit the computationual burdens of these estimators are usually heavier than those of OLS.

This paper is organized as follows. In Section 2 main results are given. Proofs of the main results are given in Section 3 and conclusions are given in Section 4.

**2. Main results.**

2.1. *Limit distributions of ordinary least square estimators.* In considering possible noise sequences $\varepsilon$, we shall impose sufficiently nice asymptotic properties of the associated partial sum processes. Let $D$ be the space of càdlàg functions $f:[0,1] \to \mathbb{R}$ equipped with the Skorokhod topology. As a measurable structure on $D$ we consider the corresponding Borel $\sigma$-algebra $\mathcal{B}(D)$. We write $X_n \xrightarrow{d}_D X$ whenever $X_n, X$ are random variables taking values in $(D, \mathcal{B}(D))$ such that $X_n$ converges weakly to $X$ in $D$ as $n \to \infty$, that is, $\lim_{n\to\infty} Ef(X_n) \to Ef(X)$ for all bounded continuous functions $f: D \to \mathbb{R}$ (see Billingsley [2] for details).

We subsume our assumptions on $\varepsilon$ into the following definition of domain of attraction.

DEFINITION 2.1. Let $\sigma^2 \in \mathbb{R}_0^+ = \mathbb{R}^+ \cup \{0\}$ and $a = (a_n)_{n \in \mathbb{N}} \subseteq \mathbb{R}^+$.

(i) We write $(\varepsilon, Z) \in \mathrm{DA}(a)$ whenever $\varepsilon = (\varepsilon_k)_{k \in \mathbb{N}}$ is a sequence of random variables and $Z$ is a random variable taking values in $(D, \mathcal{B}(D))$, satisfying the additional properties:

(D1) $$\frac{1}{a_n} \sum_{k=1}^{[n\cdot]} \varepsilon_k \xrightarrow{d}_D Z, \qquad n \to \infty.$$

(D2) $$P(Z(t) = 0 \ \forall t \in [0,1)) = 0.$$

(D3) $$\sum_{k=1}^{n} \varepsilon_k^2 = o_P(a_n^2), \qquad n \to \infty.$$

(ii) We write $(\varepsilon, Z) \in \mathrm{DA}'(a, \sigma^2)$ whenever $\varepsilon = (\varepsilon_k)_{k \in \mathbb{N}}$ is a sequence of random variables and $Z$ is a random variable taking values in $(D, \mathcal{B}(D))$, satisfying (D1) and (D2) in (i) and, in addition,

(D3)' $$\frac{1}{n} \sum_{k=1}^{n} \varepsilon_k^2 = \sigma^2 + o_P(1), \qquad n \to \infty.$$

To state a result on ordinary least squares procedures, let $Z$ a random variable taking values in $(D, \mathcal{B}(D))$. Additionally, if $Z$ satisfies condition (D2) of Definition 2.1, then a random vector $\Theta(Z)$ is well defined where

(2.6) $$\Theta(Z) = \begin{pmatrix} \left(\int_0^1 Z_s^2 \, ds\right)^{-1/2} \\ \left(\int_0^1 Z_s^2 \, ds\right)^{-1} \end{pmatrix}.$$

We define the Ornstein–Uhlenbeck process $Z_\gamma = (Z_{\gamma,t})_{t \in [0,1]}$ driven by $Z$ by setting

$$(2.7) \qquad Z_{\gamma,t} = Z_t - \gamma \int_0^t e^{-\gamma(t-s)} Z_s \, ds, \qquad t \in [0,1], \gamma \in \mathbb{R}.$$

Clearly, $Z_\gamma$ is a random variable taking values in $(D, \mathcal{B}(D))$ if $Z$ has this property. Additionally, note that (D2) is satisfied for random variables $Z$ taking values in $(D, \mathcal{B}(D))$ if and only if, for all $\gamma \in \mathbb{R}$, the same is true for the corresponding $Z_\gamma$. In this case, $\Theta(Z_\gamma)$ is well defined. Now we are in position to state the main result of this section (cf. Section 3.1 for its proof).

THEOREM 2.1. *Let $\sigma^2 \in \mathbb{R}_0^+$ and $a = (a_n)_{n \in \mathbb{N}} \subseteq \mathbb{R}^+$.*

*Let $\gamma \in \mathbb{R}$ and $(\beta_n)_{n \in \mathbb{N}} \subseteq \mathbb{R}$ satisfy $\lim_{n \to \infty} n(1 - \beta_n) = \gamma$. Let $(X_0^{(n)})_{n \in \mathbb{N}}$ be a sequence of random variables satisfying $X_0^{(n)} = o_P(a_n)$ for $n \to \infty$. As defined in (1.3), for $n \geq 1$, let $X^{(n)}$ be the corresponding nearly unstable model associated with $\beta_n$ and $\varepsilon$ with initial condition $X_0^{(n)}$.*

*Then the following assertions hold:*

(i) *If $\lim_{n \to \infty} n^{-1/2} a_n = \infty$ and $(\varepsilon, Z) \in \mathrm{DA}(a)$, then for $n \to \infty$*

$$\begin{pmatrix} a_n^{-1} n^{1/2} \hat{\tau}_n \\ n(\hat{b}_n - \beta_n) \end{pmatrix} \xrightarrow{d} \left( \frac{Z_{\gamma,1}^2}{2} + \gamma \int_0^1 Z_{\gamma,s}^2 \, ds \right) \Theta(Z_\gamma).$$

(ii) *If $\lim_{n \to \infty} n^{-1/2} a_n = 1$ and $(\varepsilon, Z) \in \mathrm{DA}'(a, \sigma^2)$, then for $n \to \infty$*

$$\begin{pmatrix} \hat{\tau}_n \\ n(\hat{b}_n - \beta_n) \end{pmatrix} \xrightarrow{d} \left( \frac{Z_{\gamma,1}^2 - \sigma^2}{2} + \gamma \int_0^1 Z_{\gamma,s}^2 \, ds \right) \Theta(Z_\gamma).$$

(iii) *If $\lim_{n \to \infty} n^{-1/2} a_n = 0$ and $(\varepsilon, Z) \in \mathrm{DA}'(a, \sigma^2)$, then for $n \to \infty$*

$$\begin{pmatrix} a_n n^{-1/2} \hat{\tau}_n \\ a_n^2 (\hat{b}_n - \beta_n) \end{pmatrix} \xrightarrow{d} -\frac{\sigma^2}{2} \Theta(Z_\gamma).$$

REMARK 2.1. (i) Note that $Z_0 = Z$. For sequences $(\beta_n) \subseteq \mathbb{R}$ with $\lim_{n \to \infty} n(1 - \beta_n) = 0$ the nearly unstable model behaves asymptotically like the strictly unstable one ($\beta_n = 1$). We may simplify the limit distributions in (i)–(iii). For instance, for $\lim_{n \to \infty} n^{-1/2} a_n = \infty$ and $(\varepsilon, Z) \in \mathrm{DA}'((a_n)_n)$ we obtain that

$$\begin{pmatrix} a_n^{-1} n^{1/2} \hat{\tau}_n \\ n(\hat{b}_n - \beta_n) \end{pmatrix} \xrightarrow{d} \tfrac{1}{2} Z_1^2 \Theta(Z).$$

(ii) For comparison of power functions of unit root tests, it is important to derive the asymptotic limit distributions under local alternatives of the form $\beta_n = 1 - \gamma/n$ (Tanaka [27], Chapter 9). In particular, Theorem 2.1 offers a first step in this direction.

REMARK 2.2. Let $\varepsilon = (\varepsilon_k)_{k \in \mathbb{N}}$ be a sequence of martingale differences satisfying (1.4) and (1.5) of Theorem 1.1. By virtue of Theorems A and B in Chan and Wei [4], we have $(\varepsilon, W) \in \mathrm{DA}'((n^{1/2})_{n \in \mathbb{N}}, 1)$, where $W$ is the Wiener process (Hall and Heyde [13], Theorem 2.23). As $W_\gamma$ is the pathwise solution of the Langevin equation, that is, almost surely,

$$dW_{\gamma,t} = -\gamma W_{\gamma,t}\, dt + dW_t, \qquad W_{\gamma,0} = 0, t \in [0,1],$$

it is easy to see by means of Itô's formula that, almost surely,

$$(2.8) \qquad \frac{W_{\gamma,1}^2 - 1}{2} + \gamma \int_0^1 W_{\gamma,s}^2\, ds = \int_0^1 W_{\gamma,u}\, dW_u.$$

Combining this identity with the fact that $(\varepsilon, W) \in \mathrm{DA}'((n^{1/2})_{n \in \mathbb{N}}, 1)$, we note that Theorem 1.1(i) is by now a consequence of Theorem 2.1(ii).

REMARK 2.3. Let $\beta_1 > 2$ and $\beta_2 > 0$ and $\varepsilon = (\varepsilon_k)_{k \in \mathbb{N}}$ be a sequence of random variables with mean zero. Additionally, suppose that $\sup_t E|\varepsilon_t|^{\beta_1 + \beta_2} < \infty$ and $L^2 = \lim_{n \to \infty} n^{-1} E(\sum_{k=1}^n \varepsilon_k)^2$ exists in $(0, \infty)$.

As shown by Phillips [20], if $\varepsilon$ is strongly mixing with mixing coefficients $\alpha_m$ satisfying $\sum_{k=1}^\infty \alpha_m^{1 - 2/\beta_1} < \infty$ (cf. Hall and Heyde [13], page 147, for the definition of mixing coefficients), then we have the result that $\sigma^2 = \lim_{n \to \infty} n^{-1} \sum_{k=1}^n E\varepsilon_k^2$ exists in $(0, \infty)$ and $(\varepsilon, LW) \in \mathrm{DA}'((n^{1/2})_{n \in \mathbb{N}}, \sigma^2)$.

By the same calculation as in Remark 2.2, Theorem 1 in Phillips [20] follows from Theorem 2.1(ii), namely, for $n \to \infty$,

$$\begin{pmatrix} \hat{\tau}_n \\ n(\hat{b}_n - \beta_n) \end{pmatrix} \xrightarrow{d} \left( \frac{W_{\gamma,1}^2 - (\sigma/L)^2}{2} + \gamma \int_0^1 (W_{\gamma,s})^2\, ds \right) \mathrm{diag}(L, 1) \Theta(W_\gamma)$$

$$= \left( \frac{1 - (\sigma/L)^2}{2} + \int_0^1 W_{\gamma,s}\, dW_s \right) \mathrm{diag}(L, 1) \Theta(W_\gamma).$$

2.2. *Applications to fractional integrated noise.* In this section we apply Theorem 2.1 to fractional integrated linear filters. We will restrict our analysis to limit processes $Z$ that are multiples of a fractional Brownian motion.

Let $H \in (0, 1)$. A centered Gaussian process $B^H = (B_t^H)_{t \in \mathbb{R}}$ with almost surely locally Hölder continuous paths of any order strictly smaller than $H$ and covariance function

$$(2.9) \qquad EB_t^H B_s^H = \tfrac{1}{2}(|t|^{2H} + |s|^{2H} - |t-s|^{2H}), \qquad s, t \in \mathbb{R},$$

is called a *fractional Brownian motion* (*FBM*) with Hurst index $H$. The choice $H = 1/2$ relates to the Wiener process $W$.

For $H \in (1/2, 1)$, FBM inherits long-range dependence in its increments, formally,

$$\sum_{t=1}^{\infty} E(B_t^H - B_{t-1}^H) B_1^H = \infty.$$

We refer to Samorodnitsky and Taqqu [25] and Marinucci and Robinson [19] for further properties and discussion.

Sequences $\varepsilon$ with partial sums in the domain of attraction of FBM have been studied by Davydov [8] and Taqqu [28] and Dobrushin and Major [10]. For more modern approaches and further examples, the interested reader is referred to Davidson and de Jong [7], Davidson [6] and Wu and Min [32]. We shall restrict our discussion to fractional integrated linear filters with i.i.d. innovations.

To this end let $\xi = (\xi_k)_{k \in \mathbb{Z}}$ be an i.i.d. sequence of random variables with $E\xi_0^2 = 1$ and $E\xi_0 = 0$. Let $\alpha = (\alpha_k)_{k \in \mathbb{N}_0} \subseteq \mathbb{R}$ be a deterministic sequence satisfying

(2.10) $$\sum_{k=0}^{\infty} |\alpha_k| < \infty.$$

Let $\eta$ be the linear filter generated by $\xi$ and $\alpha$, that is, set

(2.11) $$\eta_t = \sum_{k=0}^{\infty} \alpha_k \xi_{t-k}, \qquad t \in \mathbb{Z}.$$

For $H \in (0,1)$ we define a sequence of variables $\varepsilon^H = (\varepsilon_t^H)_{t \in \mathbb{N}}$ by setting

$$\varepsilon_t^H = (I - B)^{1/2 - H} \eta_t, \qquad t \in \mathbb{Z}, H \in (0,1).$$

Here $I$ and $B$ denote the identity and the back shift operator, respectively. The interpretation of $(I - B)^{-(H - 1/2)}$ is given by the expansion

(2.12) $$\varepsilon_t^H = \sum_{j=0}^{\infty} b_j(H) \eta_{t-j}, \qquad t \in \mathbb{Z}, H \in (0,1),$$

where

$$b_j(H) = \prod_{k=1}^{j} \frac{k - 1 + H - 1/2}{k}, \qquad j \in \mathbb{Z}, b_0(H) = 1, H \in (0,1).$$

Formally, we have $(I - B)^{H - 1/2} \varepsilon_t^H = \eta_t$; thus, $\varepsilon^H$ is called *fractional integrated noise* (see, e.g., Brockwell and Davis [3] for an overview and Robinson [23] for further discussion about Type-I and Type-II fractional integrated noise; note that the parametrization $d = H - 1/2$ is more commonly used in the econometric literature).

The sequence $\varepsilon^H$ is strictly stationary. In view of (2.10) and by means of similar calculations as in Brockwell and Davis ([3], page 470), for $H \in (0,1)\setminus\{1/2\}$ and $t \to \infty$,

$$E\varepsilon_t^H \varepsilon_1^H \sim \left(\sum_{l=0}^{\infty} E\eta_l\eta_0\right) \frac{(H-1/2)\Gamma(2-2H)}{\Gamma(3/2-H)\Gamma(H+1/2)} t^{2H-2},$$

where $\Gamma(z) = \int_0^{\infty} t^{z-1} e^{-t}\, dt$, $z > 0$, denotes Euler's $\Gamma$-function and $\sum_{l=0}^{\infty} |E\eta_l\eta_0|$ is finite.

For $H \in (1/2, 1)$ and $\sum_{l=0}^{\infty} E\eta_l\eta_0 \neq 0$, fractional integrated noise inherits long-range dependence, that is,

$$\sum_{t=1}^{\infty} E\varepsilon_t^H \varepsilon_1^H = \infty, \qquad H \in (1/2, 1).$$

Next we give sufficient conditions such that $\varepsilon^H$ and a constant multiple of $B^H$ satisfy the properties in Definition 2.1 (cf. Section 3.2 for a proof).

PROPOSITION 2.1. *Let $a, b > 2$ and $H \in (0,1)$.*

*Let $\xi = (\xi_k)_{k \in \mathbb{Z}}$ be a sequence of zero mean i.i.d. random variables with variance one and $E|\xi_0|^a < \infty$.*

*Let $\alpha = (\alpha_k)_{k \in \mathbb{N}} \subseteq \mathbb{R}$ satisfy $\sum_{k=0}^{\infty} \alpha_k^2 < \infty$ and*

$$(2.13) \qquad \sum_{k=n}^{\infty} \alpha_k^2 = O(n^{-b}), \qquad n \to \infty.$$

*Then $\sum_{k=0}^{\infty} |\alpha_k|$ is finite. If in addition $\sum_{k=0}^{\infty} \alpha_k \neq 0$, then*

$$((\varepsilon_n^H)_{n \in \mathbb{N}}, L(H, \alpha) B^H) \in \mathrm{DA}'((n^H)_{n \in \mathbb{N}}, E(\varepsilon_0^H)^2),$$

*where $\eta$ and $\varepsilon^H$ are defined in (2.11) and (2.12) and the quantities $L(H, \alpha)$ and $E(\varepsilon_0^H)^2$ are given by*

$$(2.14) \qquad L(H, \alpha) = (\Gamma(2H+1) \sin(\pi H))^{-1/2} \left|\sum_{k=0}^{\infty} \alpha_k\right|,$$

$$(2.15) \qquad E(\varepsilon_0^H)^2 = \frac{1}{2\pi} \int_{-\pi}^{\pi} \left|\sum_{k=0}^{\infty} \alpha_k e^{-ik\lambda}\right|^2 |1 - e^{-i\lambda}|^{1-2H}\, d\lambda.$$

REMARK 2.4. Proposition 2.1 can be applied to the familiar ARFIMA models. Let $H \in (0,1)$. Let $\xi = (\xi_k)_{k \in \mathbb{N}}$ admit the assumptions of Proposition 2.1. For $p, q \in \mathbb{N}$, let $\phi$ and $\theta$ be polynomials of degrees $p$ and $q$, respectively. Suppose that the polynomials have no common zeros and all of the zeros have modulus strictly larger than one. Then a sequence $\alpha = (\alpha_k)_{k \geq 0}$ is well defined by the relation $\sum_{k=0}^{\infty} \alpha_k z^k = \theta(z)/\phi(z)$, $|z| \leq 1$. As the zeros of $\phi$ have modulus strictly larger than one, there exists $r > 1$ such that $\alpha_n = O(r^{-n})$ for $n \to \infty$. Thus, $\alpha$ satisfies the assumptions in Proposition 2.1.

For this choice of $\alpha$, let $\eta$ and $\varepsilon^H$ be as defined in (2.11) and (2.12), respectively. Then $\eta$ is the unique stationary solution of $\phi(B)\eta_t = \theta(B)\xi_t$, $t \in \mathbb{Z}$, that is, an ARMA$(p,q)$ model. Furthermore, $\varepsilon^H$ is an ARFIMA$(p, H - 1/2, q)$ model (cf. Brockwell and Davis [3], pages 85–87 and Section 12.4, resp.). In view of Proposition 2.1, note that

$$((\varepsilon_n^H)_{n \in \mathbb{N}}, L(H, \alpha)B^H) \in \mathrm{DA}'((n^H)_{n \in \mathbb{N}}, E(\varepsilon_0^H)^2),$$

where we have the simplifications

$$L(H, \alpha) = (\Gamma(2H+1)\sin(\pi H))^{-1/2} \left|\frac{\theta(1)}{\phi(1)}\right|,$$

$$E(\varepsilon_0^H)^2 = \frac{1}{2\pi} \int_{-\pi}^{\pi} \left|\frac{\theta(e^{-ik\lambda})}{\phi(e^{-ik\lambda})}\right|^2 |1 - e^{-i\lambda}|^{1-2H} d\lambda.$$

We will present different representations of the limiting distribution in terms of stochastic integrals to highlight some formal similarities to Theorem 1.1. Note that $B^H$ is not a semimartingale for $H \in (0,1)\setminus\{1/2\}$. However, we may use the following fact: if $\alpha, \beta > 0$ satisfy $\alpha + \beta > 1$ and $h, g : [0,1] \to \mathbb{R}$ are Hölder continuous functions of orders $\alpha > 0$ and $\beta > 0$, respectively, then the integral $\int_0^1 g(s)\,dh(s)$ converges in the Riemann–Stieltjes sense (Young [33]). Recall that, almost surely, paths of $B^H$ and, thus, $B_\gamma^H$ are Hölder continuous of any order strictly smaller than $H$. Consequently, for $H \in (1/2, 1)$, almost surely, $\int_0^1 B_{\gamma,s}^H(\omega)\,dB_s^H(\omega)$ converges in the Riemann–Stieltjes sense. We refer to this mode of convergence as *pathwise* convergence in the sequel.

We state a direct corollary of Theorem 2.1. In particular, it holds for the noise sequences considered in Proposition 2.1.

COROLLARY 2.1. *Let $H \in (0,1)$ and $L, \sigma > 0$. Let $\gamma \in \mathbb{R}$ and $(\beta_n)_{n \in \mathbb{N}} \subseteq \mathbb{R}$ satisfy $\lim_{n \to \infty} n(1 - \beta_n) = \gamma$.*

*Let $(\varepsilon, LB^H) \in \mathrm{DA}'((n^H)_{n \in \mathbb{N}}, \sigma^2)$. Let $(X_0^{(n)})_{n \in \mathbb{N}}$ be a sequence of random variables satisfying $X_0^{(n)} = o_P(n^H)$ for $n \to \infty$. As defined in (1.3), for $n \geq 1$, let $X^{(n)}$ be the corresponding nearly unstable model associated with $\beta_n$ and $\varepsilon$ with initial condition $X_0^{(n)}$.*

*Then the following assertions hold:*

(i) *If $H \in (1/2, 1)$, then for $n \to \infty$*

$$\begin{pmatrix} n^{1/2-H}\hat{\tau}_n \\ n(\hat{b}_n - \beta_n) \end{pmatrix} \xrightarrow{d} F(H, \gamma, L, \sigma^2) \mathrm{diag}(L, 1)\Theta(B_\gamma^H),$$

*where, almost surely, $F(H, \gamma, L, \sigma^2)$ admits the representations*

$$(2.16) \quad F(H, \gamma, L, \sigma^2) = \tfrac{1}{2}(B_{\gamma,1}^H)^2 + \gamma \int_0^1 (B_{\gamma,u}^H)^2 \, du = \int_0^1 B_{\gamma,u}^H \, dB_u^H,$$

*and, almost surely, the integral on the right-hand side converges pathwise in the Riemann–Stieltjes sense.*

(ii) *If $H = 1/2$, then for $n \to \infty$*

$$\begin{pmatrix} \hat{\tau}_n \\ n(\hat{b}_n - \beta_n) \end{pmatrix} \xrightarrow{d} F(1/2, \gamma, L, \sigma^2) \operatorname{diag}(L, 1) \Theta(W_\gamma),$$

*where, almost surely, $F(1/2, \gamma, L, \sigma^2)$ admits the representations*

$$(2.17) \quad \begin{aligned} F(1/2, \gamma, L, \sigma^2) &= \frac{W_{\gamma,1}^2 - (\sigma/L)^2}{2} + \gamma \int_0^1 W_{\gamma,s}^2 \, ds \\ &= \frac{1 - (\sigma/L)^2}{2} + \int_0^1 W_{\gamma,s} \, dW_s, \end{aligned}$$

*and the integral on the right-hand side converges in the Itô sense.*

(iii) *If $H \in (0, 1/2)$, then for $n \to \infty$*

$$(2.18) \quad \begin{pmatrix} n^{H-1/2} \hat{\tau}_n \\ n^{2H}(\hat{b}_n - \beta_n) \end{pmatrix} \xrightarrow{d} -\frac{(\sigma/L)^2}{2} \operatorname{diag}(L, 1) \Theta(B_\gamma^H).$$

PROOF. In view of Theorem 2.1 and Remarks 2.2 and 2.3, it suffices to show the second identity in (2.16). Let $\gamma \in \mathbb{R}$ and $H \in (1/2, 1)$.

First note that $B_\gamma^H$ is the unique pathwise solution of

$$B_{\gamma,s}^H = -\gamma \int_0^s B_{\gamma,u}^H \, du + B_s^H, \qquad s \in [0, 1].$$

Recall that $H > 1/2$ and, almost surely, paths of $B^H$ and $B_\gamma^H$ are Hölder continuous of any order smaller than $H$. For deterministic functions with this property, Zähle [34] provides a corresponding chain rule and a density formula for Riemann–Stieltjes integrals (Zähle [34], Theorem 4.3.1 and Theorem 4.4.2). Consequently, we obtain, almost surely, that

$$\begin{aligned} \tfrac{1}{2}(B_{\gamma,1}^H)^2 &= \int_0^1 B_{\gamma,s}^H \, dB_{\gamma,s}^H = \int_0^1 B_{\gamma,s}^H \, d\left( -\gamma \int_0^s B_{\gamma,u}^H \, du + B_s^H \right) \\ &= -\gamma \int_0^1 (B_{\gamma,s}^H)^2 \, ds + \int_0^1 B_{\gamma,s}^H \, dB_s^H, \end{aligned}$$

which yields (2.16). □

REMARK 2.5. Let $H \in (1/2, 1)$, $L, \sigma^2$ and $\varepsilon = \varepsilon^H$ be the fractional integrated noise defined in (2.12). Suppose that $(\varepsilon^H, B^H) \in \mathrm{DA}'((n^H)_n, L_H, \alpha_H^2)$.

Choosing $\beta_n = 1$, $n \in \mathbb{N}$, in Corollary 2.1(i), we may additionally make use of Proposition 2.1 to relax the i.i.d. assumptions on $\eta$ as imposed in Sowell [26]. As a generalization of his result, we obtain

$$n(\hat{b}_n - 1) \xrightarrow{d} \frac{1}{2} \frac{(B_1^H)^2}{\int_0^1 (B_s^H)^2 \, ds}, \qquad n \to \infty, H \in (1/2, 1),$$

where, by now, the underlying process $\eta$ is a linear filter satisfying the assumptions of Proposition 2.1. Similar to Remark 2.1(ii), the limit in the last display remains the same for nearly unstable processes $X^{(n)}$ whenever $\beta_n$ satisfies $\lim_{n \to \infty} n(1 - \beta_n) = 0$. A related remark applies to the results of Wu [31].

REMARK 2.6. *Stationarity versus nonstationarity.* As a possible application, consider the following situation of Remark 2.4 again:

$$(1 - \beta_n B) X_t^{(n)} = \varepsilon_t^H = (1 - B)^{1/2 - H} \frac{\theta(B)}{\phi(B)} \xi_t, \qquad t \geq 1.$$

Consequently, $X^{(n)}$ is a nearly unstable ARFIMA$(p+1, H-1/2, q)$ model. For $\beta_n = 1$, the results in Corollary 2.1 can be used to test for the presence of a unit root in the autoregressive part of ARFIMA models. For $\beta_n = 1 - \gamma/n$, Corollary 2.1 gives the limit distributions under local alternatives. Although ordinary least squares estimators are easy to compute, Corollary 2.1 shows that their limit distributions are rather complicated. For $H = 1/2$, quantiles can be found in Tanaka ([27], Section 7.5). For $H \neq 1/2$, corresponding quantiles would still have to be found by simulation, which is beyond the scope of this remark.

We conclude this section with an interesting link to Malliavin's calculus. We shall restrict our discussions to fractional Gaussian noise. To this end, for $H \in (0,1)$, consider the noise sequence $\varepsilon^H = (\varepsilon_t^H)_{t \in \mathbb{N}}$ defined by $\varepsilon_t^H = B_t^H - B_{t-1}^H$, $t \in \mathbb{N}$.

Let $\gamma \in \mathbb{R}$ and $(\beta_n) \subseteq \mathbb{R}$ satisfy $\lim_{n \to \infty} n(1 - \beta_n) = \gamma$. Extending the proof of Theorem 2.1 (cf. Section 3.1), it can be seen that, for $n \to \infty$:

$$n^{-(2H \vee 1)} E \sum_{t=0}^{n-1} X_t^{(n)} \varepsilon_{t+1}^H \begin{cases} \to EF(H, \gamma, 1, 1), & H \in (1/2, 1), \, n \to \infty, \\ = 0, & H = 1/2, \, n \in \mathbb{N}, \\ \to -1/2, & H \in (0, 1/2), \, n \to \infty, \end{cases}$$

where $F(H, \gamma, 1, 1)$ is the quantity in Corollary 2.1(i).

We compute the expectation $EF(H, \gamma, 1, 1)$ for $H \in (1/2, 1)$. In Corollary 2.1(i) we have found that $F(H, \gamma, 1, 1) = \int_0^1 B_{\gamma,u}^H \, dB_u^H$, where the integral on the right-hand side converges pathwise in the Riemann–Stieltjes sense.

For $H \in (1/2, 1)$, Duncan, Hu and Pasik-Duncan [12] have developed a stochastic integral with respect to $B^H$ in the Skorokhod–Wick sense which we denote by $\int_0^1 B_{\gamma,u}^H \, \delta B_u^H$ in the sequel. Similar to the Itô integral, this integral has expectation zero. It turns out that the Riemann–Stieltjes integral in (2.16) is related to the Skorokhod–Wick integral via the relationship

$$\int_0^1 B_{\gamma,u}^H \, dB_u^H = \int_0^1 B_{\gamma,u}^H \, \delta B_u^H + \int_0^1 D_u^\phi B_{\gamma,u}^H \, du, \qquad H \in (1/2, 1), \gamma \in \mathbb{R},$$

where $D_s^\phi B_{\gamma,s}^H$ is a fractional version of the Malliavin derivative (Duncan, Hu and Pasik-Duncan [12], Theorem 3.12).

It can be easily shown that the Malliavin derivative is deterministic and, for $s \in [0, 1]$,

$$D_s^\phi B_{\gamma,s}^H = H(2H - 1) \int_0^s e^{-\gamma u} u^{2H-2} \, du, \qquad H \in (1/2, 1), \gamma \in \mathbb{R}.$$

Consequently, for $H \in (1/2, 1)$ and $\gamma \in \mathbb{R}$,

$$(2.19) \qquad EF(H, \gamma, 1, 1) = H(2H - 1) \int_0^1 \int_0^s e^{-\gamma u} u^{2H-2} \, du \, ds.$$

Particularly note that $EF(H, \gamma, 1, 1) \neq 0$ for all $H \in (1/2, 1)$ and $\gamma \in \mathbb{R}$. The integral on the right-hand side of (2.19) can be expressed in terms of incomplete $\Gamma$-functions. Direct calculations show that, for $H \in (1/2, 1)$,

$$EF(H, \gamma, 1, 1) \sim \frac{\Gamma(2H + 1)}{2} \gamma^{1-2H}, \qquad \gamma \to \infty,$$

$$EF(H, \gamma, 1, 1) \sim H(2H - 1)|\gamma|^{-2} e^{|\gamma|}, \qquad \gamma \to -\infty,$$

(2.20)
$$\lim_{H \downarrow 1/2} EF(H, \gamma, 1, 1) = \tfrac{1}{2}, \qquad \gamma \in \mathbb{R},$$

$$EF(H, 0, 1, 1) = \tfrac{1}{2}.$$

Particularly, the limit in (2.20) shows that the mapping $H \mapsto EF(H, \gamma, 1, 1)$ is not continuous at the point $1/2$. It would be interesting to extend this approach to the bias of $\hat{b}_n$ in the situation of Proposition 2.1 and compare the results to known results in the short memory case (Le Breton and Pham [17] and Larsson [16]). However, this requires substantial reasoning and will be pursued in future work.

### 3. Proofs.

3.1. *Proof of Theorem* 2.1. First we show that the corresponding nearly unstable process converges weakly to the Ornstein–Uhlenbeck process $Z_\gamma$ driven by $Z$ as defined in (2.7), that is,

$$(3.21) \qquad \frac{1}{a_n} X_{[n\cdot]}^{(n)} \xrightarrow{d}_D Z_\gamma, \qquad n \to \infty.$$

We refer to Billingsley's monograph [2] for facts on the Skorohod topology. Without loss of generality we may assume that $\beta_n > 0$ for all $n \in \mathbb{N}$ throughout the proof.

We introduce some auxiliary functions. For $n \in \mathbb{N}$, define $g_n \colon [0,1] \to \mathbb{R}$ by setting $g_n(s) = \beta_n^{-ns}$, $s \in [0,1]$, $n \in \mathbb{N}$. With the help of $g_n$, define functions $h_n \colon [0,1]^2 \to \mathbb{R}$ by

$$h_n(s,v) = \beta_n (\log \beta_n) \frac{g_n(v)}{g_n([sn]/n)}, \qquad 0 \leq s, v \leq 1, n \in \mathbb{N}.$$

Recall that $\lim_{n \to \infty} n(1 - \beta_n) = \gamma$, which implies that $\lim_{n \to \infty} n \log \beta_n = -\gamma$. As the exponential function is uniformly continuous on compacts, note that

$$(3.22) \qquad \lim_{n \to \infty} g_n(s) = \exp(\gamma s), \qquad \lim_{n \to \infty} h_n(s,u) = -\gamma e^{-\gamma(s-u)}$$

uniformly in $s, u \in [0,1]$.

For $n \in \mathbb{N}$ and $0 \leq s \leq 1$, let $Z_s^{(n)} = \frac{1}{a_n} \sum_{k=1}^{[ns]} \varepsilon_k$ be the normalized partial sum process associated with $\varepsilon$, where $a = (a_n)_{n \in \mathbb{N}}$ is the sequence as given in (D1). Let $Z_\gamma^{(n)}$ be the Ornstein–Uhlenbeck process driven by $Z^{(n)}$, as defined in (2.7). Observe that, for all $n \in \mathbb{N}$ and all $0 \leq l \leq n-1$,

$$(3.23) \qquad -n(\log \beta_n) \int_{l/n}^{(l+1)/n} Z_u^{(n)} g_n(u) \, du = \frac{\varepsilon_l}{a_n} (\beta_n^{-(l+1)} - \beta_n^{-l}).$$

On the other hand, for all $n \in \mathbb{N}$, $0 \leq t \leq n$, note that

$$(3.24) \qquad X_t^{(n)} = \beta_n^t X_0^{(n)} + \sum_{k=1}^{t} \beta_n^{t-k} \varepsilon_k.$$

In view of (3.24) and (3.23), adding and subtracting terms, the summation by parts formula yields the identity

$$(3.25) \qquad \frac{1}{a_n} X_{[ns]}^{(n)} = Z_{\gamma,s}^{(n)} + R_s^{(n)}, \qquad 0 \leq s \leq 1, n \in \mathbb{N},$$

where, for $s \in [0,1]$ and $n \in \mathbb{N}$, we set

$$R_s^{(n)} = \int_0^s (\gamma e^{-\gamma(s-u)} + h_n(s,u)) Z_u^{(n)} \, du$$
$$+ \int_{[ns]/n}^s h_n(s,u) Z_u^{(n)} \, du + \frac{X_0^{(n)}}{a_n g_n([ns]/n)}.$$

By means of (D1), observe that $\sup_{0 \leq s \leq 1} |Z_s^{(n)}|$ is bounded in probability for $n \to \infty$. Recall that $X_0^{(n)}/a_n = o_P(1)$ for $n \to \infty$. As the limits in (3.22) are uniform, we must have that $\sup_{0 \leq s \leq 1} |R_s^{(n)}| = o_P(1)$ for $n \to \infty$. In view of the continuous mapping theorem, we have $Z_\gamma^{(n)} \xrightarrow{d}_D Z_\gamma$ for $\gamma \to \infty$, giving (3.21).

Now we are in the position to prove Theorem 2.1. Note that the following identities hold, that is,

$$(3.26) \qquad \hat{\tau}_n = \frac{\sum_{t=0}^{n-1} X_t^{(n)} \varepsilon_{t+1}}{\sqrt{\sum_{t=0}^{n-1} (X_t^{(n)})^2}}, \qquad \hat{b}_n - \beta_n = \frac{\sum_{t=0}^{n-1} X_t^{(n)} \varepsilon_{t+1}}{\sum_{t=0}^{n-1} (X_t^{(n)})^2}.$$

Squaring and summing (1.3), we decompose the numerators in (3.26) into three terms which we analyze separately in the sequel, that is, for all $n \in \mathbb{N}$,

$$(3.27) \quad \sum_{t=0}^{n-1} X_t^{(n)} \varepsilon_{t+1} = \frac{1-\beta_n^2}{2\beta_n} \sum_{t=0}^{n-1} (X_t^{(n)})^2 + \frac{1}{2\beta_n}(X_n^{(n)})^2 - \frac{1}{2\beta_n} \sum_{t=1}^{n} \varepsilon_t^2.$$

We define the auxiliary random variables

$$T_{1,n} = \frac{1}{a_n^2} \sum_{t=0}^{n-1} X_t^{(n)} \varepsilon_{t+1} = \frac{n(1-\beta_n^2)}{2\beta_n} T_{2,n} + \frac{1}{2\beta_n}(T_{3,n})^2 - \frac{1}{2\beta_n} T_{4,n},$$

say, for $n \in \mathbb{N}$, where we set

$$T_{2,n} = \frac{1}{na_n^2} \sum_{t=0}^{n-1} (X_t^{(n)})^2, \qquad T_{3,n} = \frac{1}{a_n} X_n^{(n)}, \qquad T_{4,n} = \frac{1}{a_n^2} \sum_{t=1}^{n} \varepsilon_t^2.$$

As the function $f \mapsto (f(1), \int_0^1 f^2\, ds)$ is continuous as a mapping from $D[0,1]$ into $\mathbb{R}^2$, the continuous mapping theorem applies to the limit in (3.21), that is, for $n \to \infty$,

$$(3.28) \qquad \begin{pmatrix} T_{2,n} \\ T_{3,n} \end{pmatrix} = \begin{pmatrix} \int_0^1 (X_{[ns]}^{(n)}/a_n)^2\, ds \\ X_n^{(n)}/a_n \end{pmatrix} \xrightarrow{d} \begin{pmatrix} \int_0^1 Z_{\gamma,u}^2\, du \\ Z_{\gamma,1} \end{pmatrix}.$$

PROOF OF THEOREM 2.1(i). Let $\lim_{n\to\infty} n^{-1/2} a_n = \infty$ and $(\varepsilon, Z) \in \mathrm{DA}(a)$. As $\varepsilon$ satisfies property (D3) in Definition 2.1, this implies $T_{4,n} = o_P(1)$ for $n \to \infty$. Thus, by definition,

$$(3.29) \quad T_{1,n} = \frac{n(1-\beta_n^2)}{2\beta_n} T_{2,n} + \frac{1}{2\beta_n}(T_{3,n})^2 + o_P(1), \qquad n \to \infty.$$

Recall that $\beta_n \to 1$ such that $n(1-\beta_n) \to \gamma$ for $n \to \infty$. Taking this into account, we may plug the weak limit of (3.28) into (3.29) to obtain

$$(3.30) \qquad \begin{pmatrix} T_{1,n} \\ T_{2,n} \end{pmatrix} \xrightarrow{d} \begin{pmatrix} \gamma \int_0^1 Z_{\gamma,u}^2\, du + \frac{1}{2} Z_{\gamma,1}^2 \\ \int_0^1 Z_{\gamma,u}^2\, du \end{pmatrix}, \qquad n \to \infty.$$

As $Z$ satisfies (D2) in Definition 2.1, the same is true for $Z_\gamma$ for all $\gamma \in \mathbb{R}$. In particular, we obtain $\int_0^1 Z_{\gamma,s}^2 \, ds > 0$, almost surely. Combining (3.26), (3.27) and (3.30), the continuous mapping theorem implies

$$\begin{pmatrix} (n^{1/2}/a_n)\hat{\tau}_n \\ n(\hat{b}_n - \beta_n) \end{pmatrix} = T_{1,n} \begin{pmatrix} T_{2,n}^{-1/2} \\ T_{2,n}^{-1} \end{pmatrix}$$

$$\xrightarrow{d} \left( \gamma \int_0^1 Z_{\gamma,u}^2 \, du + \tfrac{1}{2} Z_{\gamma,1}^2 \right) \begin{pmatrix} (\int_0^1 Z_{\gamma,s}^2 \, ds)^{-1/2} \\ (\int_0^1 Z_{\gamma,s}^2 \, ds)^{-1} \end{pmatrix}$$

$$= \left( \gamma \int_0^1 Z_{\gamma,u}^2 \, du + \tfrac{1}{2} Z_{\gamma,1}^2 \right) \Theta(Z_\gamma), \qquad n \to \infty,$$

where $\Theta(Z_\gamma)$ is the vector as defined in (2.6), giving the result. □

PROOF OF THEOREM 2.1(ii). Let $\lim_{n\to\infty} n^{-1/2} a_n = 1$ and $(\varepsilon, Z) \in \mathrm{DA}'(a, \sigma^2)$ for some $\sigma^2 \geq 0$. We note that $T_{n,4} = \sigma^2 + o_P(1)$ for $n \to \infty$ by property (D3)$'$. Thus,

$$T_{1,n} = \frac{n(1-\beta_n^2)}{2\beta_n} T_{2,n} + \frac{(T_{3,n})^2 - \sigma^2}{2\beta_n} + o_P(1), \qquad n \to \infty.$$

In view of the limit in (3.28), the remaining proof of (ii) is similar to that of (i) and is therefore omitted. □

PROOF OF THEOREM 2.1(iii). Let $\lim_{n\to\infty} n^{-1/2} a_n = 0$ and $(\varepsilon, Z) \in \mathrm{DA}'(a, \sigma^2)$ for some $\sigma^2 \geq 0$. Note that $(a_n^2/n) T_{n,4} = \sigma^2 + o_P(1)$ for $n \to \infty$ by property (D3)$'$. Consequently, by the limit in (3.28),

$$\frac{a_n^2}{n} T_{1,n} = -\frac{1}{2}\sigma^2 + o_P(1), \qquad n \to \infty.$$

Applying analogous arguments as in (i), the proof of (iii) is complete. □

3.2. *Proof of Proposition* 2.1. If $\sum_{k=0}^\infty \alpha_k^2$ is finite and, for some $b > 2$, (2.13) holds, then $\sum_{k=0}^\infty |\alpha_k|$ is finite (cf. Wu [31] for an argument). In particularly, $\varepsilon^H$ is well defined, stationary and ergodic.

Let $\sum_{k=0}^\infty \alpha_k \neq 0$. To prove this proposition, we need to establish conditions (D1), (D2) and (D3)$'$ in Definition 2.1. For the choice $\sigma^2 = E(\varepsilon_0^H)^2$, property (D3)$'$ follows from ergodicity (cf. Brockwell and Davis [3], formulas (4.4.3) and (12.4.7) for the identity as given in (2.15)).

Property (D2) is obvious. In order to show that (D1) holds, we shall make use of Theorem 3.1 of Davidson and de Jong [7].

Therefore it suffices to establish their conditions (b)–(d). In view of Jensen's inequality, as $E|\xi|^b$ is finite, note that $\sup_t E|\eta_t|^b$ is finite, giving (b) in

Davidson and de Jong [7]. For all $m \in \mathbb{N}$ and $t \in \mathbb{Z}$, observe that

$$E(\eta_t - E[\eta_t|\xi_{-m},\ldots,\xi_0])^2 \leq \sum_{k=m+1}^{\infty} \alpha_k^2.$$

In view of the assumption (2.13), the inequality in the last display implies that $\eta$ is $L_2$-NED on $\xi$ of size $-1/2$, giving (c) in Davidson and de Jong [7] (cf. de Jong and Davidson [9] for the definition of NED). In view of the dominated convergence theorem, it is straightforward to show

$$\lim_{n\to\infty} \frac{1}{n} \sum_{t=1}^{n} \sum_{s=1}^{n} E(\eta_t \eta_s) = \left(\sum_{k=0}^{\infty} \alpha_k\right)^2 \in (0,1),$$

giving (d) in Davidson and de Jong [7]. Consequently, Theorem 3.1 of Davidson and de Jong [7] applies, that is,

$$\left(E\left(\sum_{k=1}^{n} \varepsilon_k^H\right)^2\right)^{-1/2} \sum_{k=1}^{[n\cdot]} \varepsilon_k^H \xrightarrow{d}_D B^H, \qquad n \to \infty.$$

Combining Lemma 3.1 of Davidson and de Jong [7] and equation (9.3) in Doukhan, Oppenheim and Taqqu ([11], page 28), this yields

$$E\left(\sum_{k=1}^{n} \varepsilon_k^H\right)^2 \sim L(H,\alpha)^2 n^{2H}, \qquad n \to \infty,$$

where $L(H,\alpha)$ is the quantity given in (2.14), giving the result.

**4. Conclusion.** In this paper, asymptotic distributions of the least squares estimator of a nearly unstable AR(1) process with long-memory errors are derived. It is shown that the limiting distributions behave very differently from those in the short-range dependent case, although formal similarities can be established. Nevertheless, the results shed light on the trickiness of dealing with long-memory models. With these results, we can gather the necessary framework to tackle the more challenging problem of a general nearly unstable AR($p$) model in the spirit of Jeganathan [14]. This will be pursued in a forthcoming paper. Finally, the asymptotic theory developed here can be used to study the important issue of fractionally cointegrated systems in a nearly unstable environment.

**Acknowledgments.** The authors takes pleasure in thanking the Editor, an Associate Editor and two anonymous referees for their insightful suggestions and for providing us with relevant literature. This research was supported in part by HKSAR-RGC Grants CUHK400305 and CUHK400306. Part of this research was completed when the first author visited the Institute of Mathematical Sciences (IMS) at CUHK in 2004 and 2006. Research support from the IMS is gratefully acknowledged.

CMA AND MASCOS  
AUSTRALIAN NATIONAL UNIVERSITY  
CANBERRA, ACT 0200  
AUSTRALIA  
E-MAIL: Buchmann@maths.anu.edu.au

DEPARTMENT OF STATISTICS  
CHINESE UNIVERSITY OF HONG KONG  
SHATIN, NT  
HONG KONG  
E-MAIL: nhchan@sta.cuhk.edu.hk